\documentclass[a4paper,11pt]{article}
\usepackage{geometry,amssymb,amsmath,enumerate,latexsym,array}
\usepackage[dvips]{graphicx}
\usepackage[all]{xy}
\begin{document}
\newtheorem{example}{Example}[section]
\newtheorem{rem}[example]{Remark} 
\newtheorem{lem}[example]{Lemma} 
\newtheorem{fact}[example]{Fact} 
\newtheorem{thm}[example]{Theorem} 
\newtheorem{cor}[example]{Corollary} 
\newtheorem{Def}[example]{Definition} 
\newtheorem{prop}[example]{Proposition} 
\newenvironment{proof}{\noindent {\bf Proof} }{\mbox{} \hfill $\Box$
\mbox{}\\}

\title { Gr\"obner Basis Procedures for Testing Petri Nets\\
\thanks{
AMS 1991 CLASSIFICATION: }}
\author{ Angie Chandler \thanks{Supported 1996-99 by an EPSRC Research
Studentship.}
\\Engineering Department \\ Faculty of Applied Science \\ 
Lancaster University \\ Lancaster LA1 4YR \\ United Kingdom 
\\a.k.chandler@lancaster.ac.uk
\and Anne Heyworth\thanks{ Supported 1995-8 by an EPSRC
Earmarked Research Studentship, `Identities among relations for
monoids and categories', and 1998-9 by  a University of Wales, Bangor,
Research Assistantship.}\\ School of Mathematics \\ University of Wales,
Bangor \\ Dean Street, Bangor\\
Gwynedd LL57 1UT \\ United Kingdom \\ a.l.heyworth@bangor.ac.uk
 } \maketitle
 
\begin{abstract}
This paper contains introductory material on Petri nets and
Gr\"obner basis theory and makes some observations on the relation
between the
two areas. The aim of the paper is to show how Gr\"obner basis
procedures can
be applied to the problem of reachability in Petri nets, and to give
details 
of an application to testing models of navigational systems.
\end{abstract}

\section{Introduction}

{\em Petri nets} are a graphical and mathematical modelling tool applicable to
many 
systems. They may be used for specifying information processing systems
that 
are concurrent, asynchronous, distributed, parallel, non-deterministic, 
and/or stochastic. 
Graphically, Petri nets are useful for illustrating and describing
systems, and  
tokens can simulate the dynamic and concurrent activities. 
Mathematically, it is possible to set up models such as state equations
and 
algebraic equations which govern the behaviour of systems. 
Petri nets are understood by both practitioners and theoreticians and so
provide a 
powerful communication link between them. 
For example, engineers can show mathematicians how to make practical and 
realistic models and mathematicians may be able to produce theories to
make 
the systems more methodical or efficient, which is in fact demonstrated
by this
collaborative paper.\\

The area of computer algebra 
called {\em Gr\"obner basis theory} includes the rewriting theory widely 
used in computer science and provides methods for  handling the rule 
systems defining various types of algebraic structure. 
It has been proved that it is not always possible to deduce all 
consequences of a system of rules --
when it is possible the levels of complexity involved quickly
require the use of computers.
In the commutative case computational Gr\"obner basis methods
have has been successfully applied in theorem proving, robotics,
integer programming, coding theory, signal processing, enzyme
kinetics, experimental design, differential equations, and many
others. All major computer algebra packages now include
implementations of these procedures, and pocket calculator
implementations will soon be available. A collection of recent papers on
Gr\"obner basis research is \cite{RISC}.\\

In this paper we show how Gr\"obner basis procedures can be applied to
reversible Petri nets to solve the reachability problem. This provides a 
practical test which can be useful in the design and analysis of Petri
nets.
In particular the examples show a practical application of the Gr\"obner
basis
methods to Petri nets modelling navigation systems. Further details of 
these mechatronic navigation systems can be found in \cite{Angie}.
Related algebraic research, and preliminaries to this paper may be found 
in \cite{Anne}.

\section{Background to Gr\"obner Bases}

We give a brief summary of the main results in commutative Gr\"obner
basis theory that will be used in this paper.
For a fuller introduction to the subject see \cite{AdLo,Co+}.\\

Let $X$ be a set. Then the elements of $X^\Delta$ are all
power products of elements of $X$, including an identity $1$, with
multiplication defined in the usual way. The commutativity
condition is summarised by $xy=yx$ for all $x,y \in X$. Let $K$ be a
field (the field of rational numbers, $\mathbb{Q}$ suffices for our work). 
Then $K[X^\Delta]$ is the ring of commutative
polynomials
$$
f=k_1m_1+ \cdots + k_tm_t
$$
where $k_1, \ldots ,k_t \in K$ and $m_1, \ldots, m_t \in X^\Delta$
with the operations of polynomial addition and polynomial multiplication 
defined in the usual way.\\

Consider a set of polynomials $P \subseteq K[X^\Delta]$. 
We say that two polynomials $f$ and $g$ of $K[X^\Delta]$ are
{\em equivalent modulo $P$} and write $f =_P g$ if their
difference can be expressed in terms of $P$, i.e.
$$
f-g=u_1p_1 + \cdots + u_np_n
$$ 
for some $p_1,\cdots,p_n \in P, u_1, \cdots, u_n \in K[X^\Delta]$.\\

In 1965 Bruno Buchberger invented the concept of a Gr\"obner basis \cite{Bu65}. 
Techniques of Gr\"obner basis theory enable us to decide whether or not
$f =_P g$ for given $P$, $f$, $g$ in $K[X^\Delta]$ as above.\\

Computation begins by specifying an ordering $>$ on the power products
(this must be a well-ordering, compatible with multiplication).
This enables us to define reduction modulo a set of polynomials $P$ --
multiples of polynomials in $P$ are subtracted from a given polynomial
$f$ in order to obtain successively smaller polynomials --
the reduction is denoted $\to_P$.
The reflexive, symmetric, transitive closure of $\to_P$ coincides with the
congruence $=_P$. If $P$ is a Gr\"obner basis then $\to_P$ is
confluent, meaning that there is a unique irreducible element in each
congruence class, obtainable from any other element by repeated reduction 
modulo $P$. If $P$ is not a Gr\"obner basis then it is always possible to
use Buchberger's algorithm to obtain a set of polynomials $Q$ which is a 
Gr\"obner basis such that $=_P$ coincides with $=_Q$.\\

Thus, given a set of polynomials $P \subseteq K[X^\Delta]$, the problem of 
deciding whether $f$ is equivalent to $g$ modulo $P$ for any $f, g$ in
$K[X^\Delta]$ can always be determined by calculating a Gr\"obner basis $Q$.
The polynomials are equivalent if and only if their difference $f-g$
reduces modulo $Q$ to zero.\\

We will not explain these calculations in any greater detail, 
but refer the reader to texts on Gr\"obner bases, such as \cite{AdLo,Co+}.
In the commutative case it is always possible to determine a Gr\"obner 
basis, but computers are usually required for all but the most
basic problems.
In our examples we use $\mathsf{MAPLE}$ and $\mathsf{GAP3}$, with some 
Gr\"obner basis procedures implemented by the second author \cite{Anne}.

\section{Petri Nets} 

A Petri net has two types of vertices:
\emph{places} (represented by circles) and 
\emph{transitions} (represented by double lines). Edges exist only between places
and transitions and are labelled with their weights. 
In modelling, places represent conditions and transitions represent
events. 
A transition has input and output places, which represent preconditions
and 
postconditions (respectively) of the event. 
A good introduction to the ideas of Petri nets is \cite{Murata}. 
 
\begin{Def}[Petri Net]
A \textbf{Petri net} (without specific initial marking) is a quadruple 
$\underline{N}=(X,T,\mathcal{F},w)$ where:
$X$ is a finite set (of places), 
$T$ is a finite set (of transitions), 
$\mathcal{F} \subseteq (X \times T) \cup (T \times X)$ is a set of edges
(flow relation) and
$w:\mathcal{F} \to \mathbb{N}$ is a weight function.
\end{Def}

The state of a system is represented by the assignation of ``tokens'' 
to places in the net.

\begin{Def}[Marking]
A \textbf{marking}  is a function $M:X \to \mathbb{N} \cup \{0\}$.
\end{Def}

Dynamic behaviour is represented by changes in the state of the Petri net 
which is formalised by the concept of firing.

\begin{Def}[Firing Rule] \mbox{ }\\
\begin{enumerate}[i)] 
\item 
A transition $t$ is \textbf{enabled} if each input place $x$ of $t$ is
marked with at least $w(x,t)$ tokens.
\item 
An enabled transition may or may not \textbf{fire} -- depending on
whether or not the relevant event occurs. 
\item 
Firing of an enabled transition $t$ removes $w(x,t)$ tokens from each 
input place $x$ of $t$ and adds $w(t,y)$ tokens to each output place $y$
of $t$. 
\end{enumerate}
\end{Def}

Despite their apparant simplicity, Petri nets can be used to model complex
situations -- for some examples see \cite{Desel}.
One of the main problems in Petri net theory is reachability -- the problem
corresponds to deciding which situations (modelled by the net) are possible,
given some sequence of events.

\begin{Def}[Reachability]
A marking $M_1$ is said to be \textbf{reachable} from a marking $M_2$ in
a net $\underline{N}$, 
if there is a sequence of firings that transforms $M_2$ to $M_1$.
Often a Petri net comes with a specified {\em initial marking} $M_0$. 
The \textbf{reachability problem} for a Petri net $\underline{N}$ with 
initial marking $M_0$ is:
{\em Given a marking $M$ of $\underline{N}$, 
is $M$ reachable in $\underline{N}$?}
\end{Def}
  
For the type of Petri nets defined so far, reachability is
decidable in exponential time and space \cite{Murata}. \\

Reversibility is a property of Petri nets corresponding to the potential for 
the device being modelled to be reset. For our applications it is essential 
that we can reset, therefore this property is vital.

\begin{Def}[Reversibility]
A Petri net $\underline{N}$ is called
\textbf{reversible} if a marking $M'$
is reachable from a marking $M$ in $\underline{N}$, 
then $M$ is reachable from $M'$.
\end{Def}

Different definitions of reversibility exist. The definition we use
is chosen for engineering rather than mathematical reasons as
in \cite{Murata}. The paper \cite{Caprotti} by Caprotti, 
Ferscha and Hong contains a result apparently similar to ours, but they use
a different definition of reversibility, which is much more 
restrictive -- perhaps this is appropriate for different applications.\\

In order to apply Gr\"obner basis techniques
we use monomials to represent the markings (there is a one-to-one
correspondence between monomials and markings), and so associate
a transition with the difference between two monomials (input and output).

\begin{Def}[Polynomial Associated with a Marking]
Let $\underline{N}=(X,T,\mathcal{F},w)$ be a Petri net. To every marking
$M$ we will associate a polynomial 
$$
pol(M):=\prod_X \; x^{M(x)},
$$ 
that is the formal product of elements of $X$ raised to the power $M(x)$
(the number of tokens held at the place $x$).
\end{Def}

\begin{Def}[Polynomial Associated with a Transition] 
Each transition $t$ has an associated polynomial 
$$
pol(t):= \prod_X \; x^{w(x,t)} - \prod_X \; y^{w(t,y)},
$$ 
that is the input required for the transition to be enabled
minus the output resulting from a firing. We often write $pol(t)=l-r$,
to distinguish the two terms.
\end{Def}

To represent the dynamic structure we must consider how the
transition polynomials are related to polynomials of markings which enable 
them and how firings of transitions affect the polynomials of the markings. 
Suppose a marking $M_i$ enables a transition $t_i$. By the definitions it is 
clear that this corresponds to $pol(M_i)$ being equal to $u_i \, l_i$ where
$pol(t_i)=l_i-r_i$ and $u_i$ is a power product in $X^\Delta$. 
It then follows that if $t_i$ fires, the resulting marking $M_{i+1}$ will
have polynomial $pol(M_{i+1})=pol(M_i)-u_i\,pol(t_i)=u_i \, r_i$.

\newpage
\begin{example}[Polynomials and the Firing Rule]
\mbox{ }

\begin{center}
\includegraphics[width=10cm]{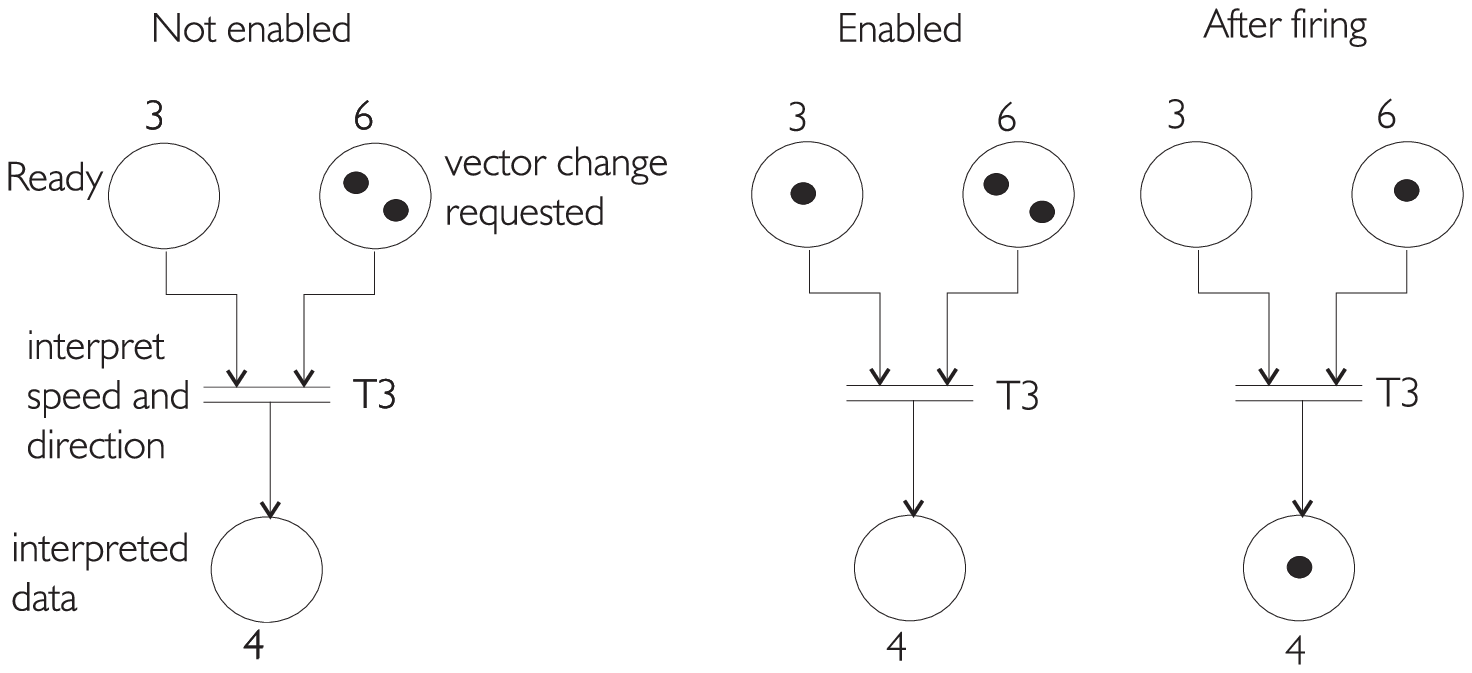}
\end{center}

\emph{The diagrams above show three different states of a transition 
$t_3$ of a Petri net Example \ref{motors}.
The polynomial associated with the transition is $pol(t_3)=x_3x_6-x_4$.
The first marking $M_1$ does not enable $t_3$; this corresponds to the
fact that $pol(M_1)=(x_6)^2$ is not a multiple of $x_3x_6$. 
The second marking $M_2$ does enable $t_3$, and $pol(M_2)=x_3(x_6)^2$.
The marking 
resulting from the firing of $t_3$ after it has been enabled by $M_2$ is
$M_3$.
In terms of polynomials the firing is represented by 
$pol(M_3)=pol(M_2)-x_6pol(t_3)=x_4x_6$.
A firing sequence is denoted by  
$M_0 \stackrel{t_1}{\to} M_1 \stackrel{t_2}{\to} \cdots
\stackrel{t_n}{\to} M_n$ where 
the $M_i$ are markings and the $t_i$ are transitions (events)
transforming 
$M_{i-1}$ into $M_i$. 
In terms of polynomials the above firing sequence gives the information
$pol(M_n) = pol(M_0) - u_1\,pol(t_1) - u_2\,pol(t_2) - \cdots - u_n\,pol(t_n)$
for some $u_1, u_2, \ldots, u_n \in X^\Delta$.}
\end{example}

\begin{thm}[Reachability and Equivalence of Polynomials]\mbox{ }\\
Let $\underline{N}$ be a reversible Petri net with initial marking
$M_0$.
Define $P:=\{pol(t) : t \in T\}$.
Then a marking $M$ is reachable in $\underline{N}$ if and only if
$pol(M_0) =_P pol(M)$.
\end{thm}

\begin{proof}
First suppose that $M$ is reachable.
Then there is a firing sequence
$M_0 \stackrel{t_1}{\to} M_1 \stackrel{t_2}{\to} \cdots 
\stackrel{t_{n-1}}{\to} M_{n-1} \stackrel{t_n}{\to} M$.
Therefore, as above, there exist $u_1, \ldots, u_n \in
X^\Delta$ such that
$pol(M_0)-pol(M) = u_1pol(t_1) + \cdots + u_n pol(t_n)$.
Hence $pol(M_0) = _P pol(M)$.\\

For the converse, suppose $pol(M_0) =_P pol(M)$.
Then 
$$
pol(M_0)  = pol(M) \pm u_1pol(t_1) \pm \cdots \pm u_mpol(t_m).
$$
The proof is by induction on $m$.\\

For the base step put $m=0$ then $pol(M_0)=pol(M)$. The correspondence between markings 
and their associated polynomials is one-to-one, so here $M_0=M$ and $M$ is
clearly reachable.\\

For the induction step we assume that a marking $M'$ is reachable from $M_0$ if
$$
pol(M_0)  = pol(M') \pm u_1pol(t_1) \pm \cdots \pm u_{m-1}pol(t_{m-1}).
$$
for a fixed $m$. Now suppose $M$ is a marking such that
$$
pol(M_0)  = pol(M) \pm u_1pol(t_1) \pm \cdots \pm u_mpol(t_m).
$$
Then for some $i \in \{1, \ldots, m\}$
either $pol(M_0)=u_il_i$ or $pol(M_0)=u_ir_i$ where $pol(t_i)=l_i-r_i$.

In the first case $pol(M_0)=u_il_i$.
Observe that $M_0$ enables $t_i$ and define a marking $M'$
by $M_0 \stackrel{t_i}{\to} M'$. Then 
$$
pol(M')=pol(M) \pm u_1pol(t_1) \pm \cdots \pm u_{i-1}pol(t_{i-1}) \pm
 u_{i+1}pol(t_{i+1}) \pm \cdots \pm u_mpol(t_m)
$$
so, by assumption, $M$ is reachable from $M'$ and so $M$ is reachable from $M_0$.

In the second case $pol(M_0)=u_ir_i$.
There is a marking $M'$ such that $pol(M')=u_ir_i$ and
$$
pol(M')=pol(M) \pm u_1pol(t_1) \pm \cdots \pm u_{i-1}pol(t_{i-1}) \pm
u_{i+1}pol(t_{i+1}) \pm \cdots \pm u_mpol(t_m).
$$
Now, $M$ is reachable from $M'$ by assumption and $M_0$ is reachable from $M'$ 
by a firing of $t_i$. By reversibility, therefore, $M'$ is reachable from $M_0$
and hence $M$ is reachable from $M_0$. 
\end{proof}

\begin{cor}[Gr\"obner Bases Determine Reachability]\mbox{ }\\
Reachability in a reversible Petri net can be determined using a
Gr\"obner basis.
\end{cor}
\begin{proof}
Let $K$ be a field.
First observe that $P \subseteq K[X^\Delta]$.
Let $Q$ be a Gr\"obner basis for $P$.
Then $pol(M)=_P pol(M_0)$ if and only if there exists $p \in
K[X^\Delta]$ such that $pol(M)$ and 
$pol(M_0)$ reduce to $p$ by $\to_Q$. 
\end{proof}

\begin{rem}[Catalogue of Reachable Markings]
\emph{Recall that Gr\"obner bases techniques use an ordering on the
power
products. There is a one-to-one correspondence between power products and
markings. We can begin to catalogue the markings in increasing order.
Given a Gr\"obner basis for the polynomials of the transitions of a
Petri net it can be determined whether each marking is reachable: if the
power product reduces to the same irreducible power product as the
initial marking then it is reachable. In this way the Gr\"obner basis
can be used to build up a list of reachable markings.}
\end{rem}

\begin{rem}[Testing for Reversibility in Petri Net Design]
\emph{The reversibility of a Petri net can be interpreted as the
ability to
reset the application it models. Whilst the reachability of a place,
given an initial marking, can be determined by standard means,
reversibility cannot be established directly.}\\

\emph{Calculating a Gr\"obner basis for the Petri net makes the
determination
of reachable markings much more obvious, and unwanted markings can be
immediately detected. There are two reasons why unwanted markings may
occur.
In the first case there is a basic error in the net which allows some
firing sequence of marking which should be avoided; the Gr\"obner basis
is effective in showing up these markings. The second type of problem
occurs when marking supposed to be unreachable is found to be reachable,
the implication here being that the net is not truly reversible. As
reversibility is a desirable property, the net can then be modified and
retested.}\\

\emph{In practical terms Gr\"obner bases have been shown by the authors
to be
useful in Petri net design -- repeated testing by computing Gr\"obner
bases shows up unintended effects or non-reversibility.
Our examples are Petri nets designed by the first author to model
software interfaces to hardware components of mobile robot navigation systems, 
and their development was helped in this way.}
\end{rem}

\begin{example}[Software Interface for Motors]\label{motors}
\emph{This Petri net represents the software interface between a
user and the set of motors used to drive a mobile robot.}

\begin{center}

\includegraphics[ height=10cm]{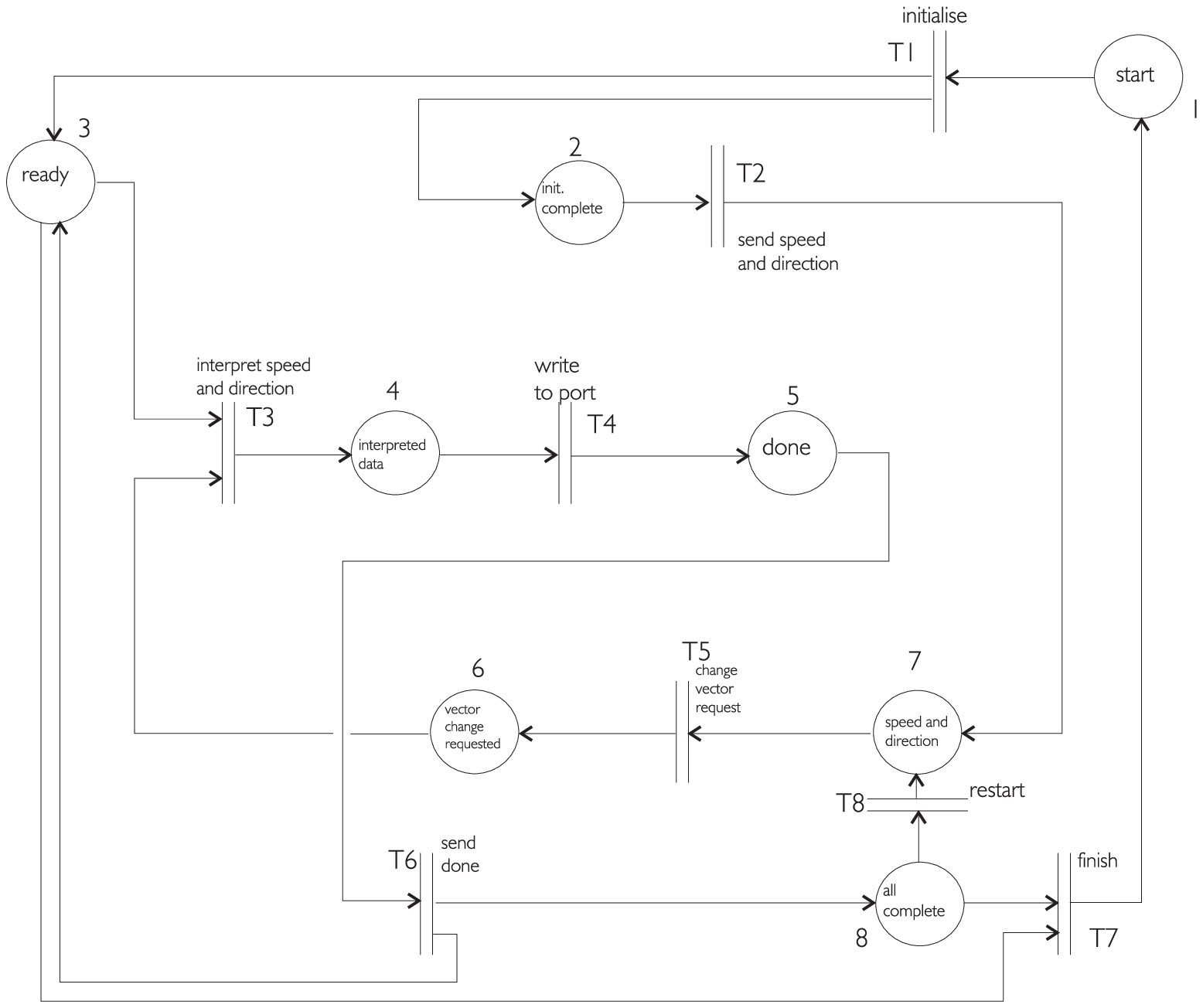}

\end{center}

\emph{Here, once the motors have been initialised, the user may input
the
required
speed and direction for each motor.  This information is then
interpreted
and written to the relevant port, if there is also a token available in 
the ``ready'' place (3), to enable the ``interpret speed and direction''
transition $t_3$.}

\emph{The places are labelled $x_1,\ldots,x_{11}$. There are eight
transitions, and their polynomials are as follows:}
\begin{center}
\begin{tabular}{llllllll}
$pol(t_1)=$ & $x_1-x_2x_3$ & $pol(t_2)=$ & $x_2-x_7$ & 
$pol(t_3)=$ & $x_3x_6-x_4$ & $pol(t_4)=$ & $x_4-x_5$ \\
$pol(t_5)=$ & $x_7-x_6$    & $pol(t_6)=$ & $x_5-x_3x_8$ &
$pol(t_7)=$ & $x_3x_8-x_1$ & $pol(t_8)=$ & $x_8-x_7$
\end{tabular}
\end{center}

{\em The Gr\"obner basis for this set of polynomials -- with respect to a 
degree-lexicographic ordering -- is} 
$$\{x_4-x_1,\ x_5-x_1,\ x_6-x_2,\ x_7-x_2,\ x_8-x_2,\ x_2x_3-x_1\}.$$

{\em The catalogue of markings reachable from an initial marking $x_1$ 
is quickly calculated to be:}
$$\{x_1,\ x_4,\ x_5,\ x_2x_3,\ x_3,\ x_6,\ x_3x_7,\ x_3x_8\}.$$
{\em This catalogue can be examined by the Petri net designer who
interprets the different states. When unexpected states appear in the
catalogue it indicates an error, which generally signifies that the net
is not reversible.}\\

\emph{For Petri nets such as this to execute efficiently, 
it is essential that the user can confirm both the reachability 
and the reversibility of the net.  For instance, should the place 
``done" (5) prove to be unreachable from an initial marking where the 
place ``start" (1) held a token, this would show that no data would be 
written to the port in transition ``write to port" ($t_4$),
thus making the motors uncontrollable.  
If the net here was non-reversible, it would indicate that the motors
could not be disabled, which in this situation is undesirable.
Once the Petri net has been tested for such bugs, the user need only 
concern themselves with the simple functions executed within individual 
transitions, greatly decreasing the likelihood of a serious, or perhaps 
dangerous, failure of the robot.}
\end{example}

\section{Coloured Petri Nets}

A coloured Petri net circulates tokens of more than one
type. The transitions in the net are affected differently by
different combinations of colours of tokens.
An example of this is where tokens represent data signals. 
Incomplete or corrupt signals should be dealt with differently from
complete signals, these two types of data would be represented by
different colours of tokens (``pass'' and ``fail'' in Example
\ref{compass}).\\ 

Recall that if $C$ is a set (of colours) then $C^\Delta$ is the set of
all power products of elements of $C$. Essentially an element of
$C^\Delta$ assigns a non-negative integer to each element of $C$. The
definition of a coloured Petri net that we give uses this kind of
notation, but is equivalent to that given by Murata in \cite{Murata}.
One element $m$ of $C^\Delta$ is said to be a {\em multiple} of another
element $l$ if
$m=ul$ for some $u \in C^\Delta$. 

\begin{Def}[Coloured Petri Net]
A \textbf{coloured Petri net} is a quintuple $\underline{N}_C=
(X,T,C,\mathcal{F},w)$, where $X$ is a set of places, $T$ is a set of
transitions, $C$ is a set of colours, $\mathcal{F} \subseteq (X \times
T) \cup (T \times X)$ is the flow relation and $w : \mathcal{F} \to
C^\Delta$. A \textbf{marking} in $\underline{N}_C$ is a function 
$M: X \to C^\Delta$.
The \textbf{firing rule} is as follows:
\begin{enumerate}[i)]
\item
A transition $t$ is enabled if each input place $x$ of $t$ is marked
with a multiple of $w(x,t)$. 
\item
An enabled transition may or may not fire.
\item
A firing of an enabled transition $t$ deletes the power product $w(x,t)$
from the marking at each input place $x$, and appends the marking at
each output place $y$ with the power product $w(t,y)$. 
\end{enumerate}
\end{Def}

A coloured Petri net can in fact be considered as a structurally folded
version of an ordinary Petri net if the number of colours is finite.
Each place $x$ is unfolded into a set of places, one for each colour of
token which $x$ may hold, and each transition $t$ is unfolded into a
number of transitions, one for each way that $t$ may fire.
It is immediate that the techniques discussed in the previous section
may be applied to coloured Petri nets. In fact we can pass directly from
the coloured Petri net to commutative polynomials in 
$K[(X \times C)^\Delta]$, where $K$ is a field. 
Elements of $(X\times C)^\Delta$ are written 
$(x_1,c_1)\cdots(x_n,c_n)$, where $x_1,\ldots,x_n \in X$, and 
$c_1, \ldots, c_n \in C$. We define $(x_i,c_i)(x_j,c_j)=(x_i,c_ic_j)$
when
$x_i=x_j$.\\

\begin{thm}[Gr\"obner Bases for Coloured Petri Nets]\mbox{ }\\
Let $\underline{N}_C$ be a coloured Petri net.
If $M$ is a marking in $\underline{N}_C$, then define
the {\em polynomial associated with the coloured marking} to be
$pol(M):=\prod_X (x,M(x))$.
Similarly if $t$ is a transition in $\underline{N}_C$, then define
the {\em polynomial associated with the coloured transition} to be
$pol(t):=\prod_X (x,w(x,t))- \prod_X (y,w(t,y))$.\\

From these definitions we observe that a transition $t$ in a coloured 
Petri net has an associated polynomial of the form
$pol(t)=l-r$ where $l,r \in (X \times C)^\Delta$. The transition $t$
is enabled by a marking $M$ if $pol(M)=ul$, for some 
$u \in (X \times C)^\Delta$.
If $t$ fires then the new marking has associated polynomial $ur$.

It follows that if we
define $P:=\{pol(t):t \in T\}$ then a marking $M$ is reachable 
if and only if $pol(M) =_P pol(M_0)$.
Therefore if $Q$ is a Gr\"obner basis for $P$ it is decidable whether or not
$M$ is reachable in $\underline{N}$.
\end{thm}

The results (and proofs) are naturally very similar to the results
for standard Petri
nets. The value is in the application -- where it is more efficient to
work with coloured nets it is appropriate to associate polynomials to these
models directly.

\begin{example}[Software Interface for Compass]\label{compass}
\emph{The following Petri net shows the software interface to an
external 
compass, where the compass provides data in the form of 
an ASCII string.}

\emph{The states here are numbered, but two types of token: ``pass'' ($x$) 
and ``fail'' ($y$), circulate in the net. 
This Petri net is initialised with a single 
``pass" ($x$) token at the ``start" place ($1$) 
together with a ``pass" ($x$) and a ``fail" ($y$) token in each of the places
``input" ($18$) and ``continue" (${19}$).  
The additional tokens at (${18}$) and (${19}$)
provide the colouring    
essential for rigorous testing of this Petri net.  For instance, when the
``return data" $t_3$ or $t_{20}$ transition is fired, the colour of the
token  
output to place ``raw data ready" ($3$) depends solely on the colour of
the token from place ``input" (${18}$).  The transitions ``read in"
$t_4$ or
$t_{21}$, ``calculate checksum" $t_5$ or $t_{22}$ and ``test" $t_6$ or
$t_{23}$ 
will output a token matching the input token, having no effect on the
colouring, but
the transition ``find bearing" $t_7$ will only be enabled by a ``pass"
token,
which represents a received ASCII string with a correct checksum, as    
determined in the ``test" $t_6$ or $t_{23}$ transition.  A ``fail" token
would  
instead enable the transition ``data request" $t_{16}$, which will
provide a
value using dead reckoning in place of the corrupted data.}\\

\emph{Colouring of this net is helpful, as it ensures that only complete
uncorrupt data is used. The Petri net of this example was constructed by
repeated testing using Gr\"obner basis methods.
We use $x_i$ to denote a ``pass'' token at place $i$, and $y_i$ to denote a 
``fail'' token at place $i$. The initial marking is therefore associated with 
the monomial $x_1x_{18}x_{19}y_{18}y_{19}$.
The set $P$ of polynomials associated with the transitions is as follows:}
\begin{center}
$x_1-x_2x_4,~ x_5-x_{12},~ x_2x_{18}-x_3x_{18},~ y_2y_{18}-y_3y_{18},~
x_3x_{13}-x_6,~ y_3x_{13}-y_6, x_6-x_7, y_6-y_7,$\\
$ x_7-x_8,~ y_7-y_8,~
x_8-x_{10},~ x_{12}-x_{13},~ x_{11}-x_2x_{14},~ x_{14}x_{19}-x_{15}x_{19},~
x_{14}y_{19}-y_{15}y_{19},~ y_{15}-x_{17},$\\
$x_3x_{17}-x_{16},~ y_3x_{17}-x_{16},~ x_2x_{17}-x_{16},~ x_{15}-x_{12},~ 
x_4-x_5,~ y_8-x_9,~ x_9-x_{11},~ x_{10}-x_{11},~ 
x_{16}-x_1.
$
\end{center}

\emph{Using $\mathsf{MAPLE}$ a Gr\"obner basis $Q$ for $P$ with respect to 
the order $tdeg$ has 47 rules}

\emph{Given the initial marking $x_1x_{18}x_{19}y_{18}y_{19}$, 
there are 11 reachable markings having five tokens and 32 reachable markings having
six tokens. Examining the catalogue of reachable states and relating them to 
the situations they represent will confirm that the net will behave as the user
would expect.}

\begin{center}

\includegraphics[ height=21.5cm]{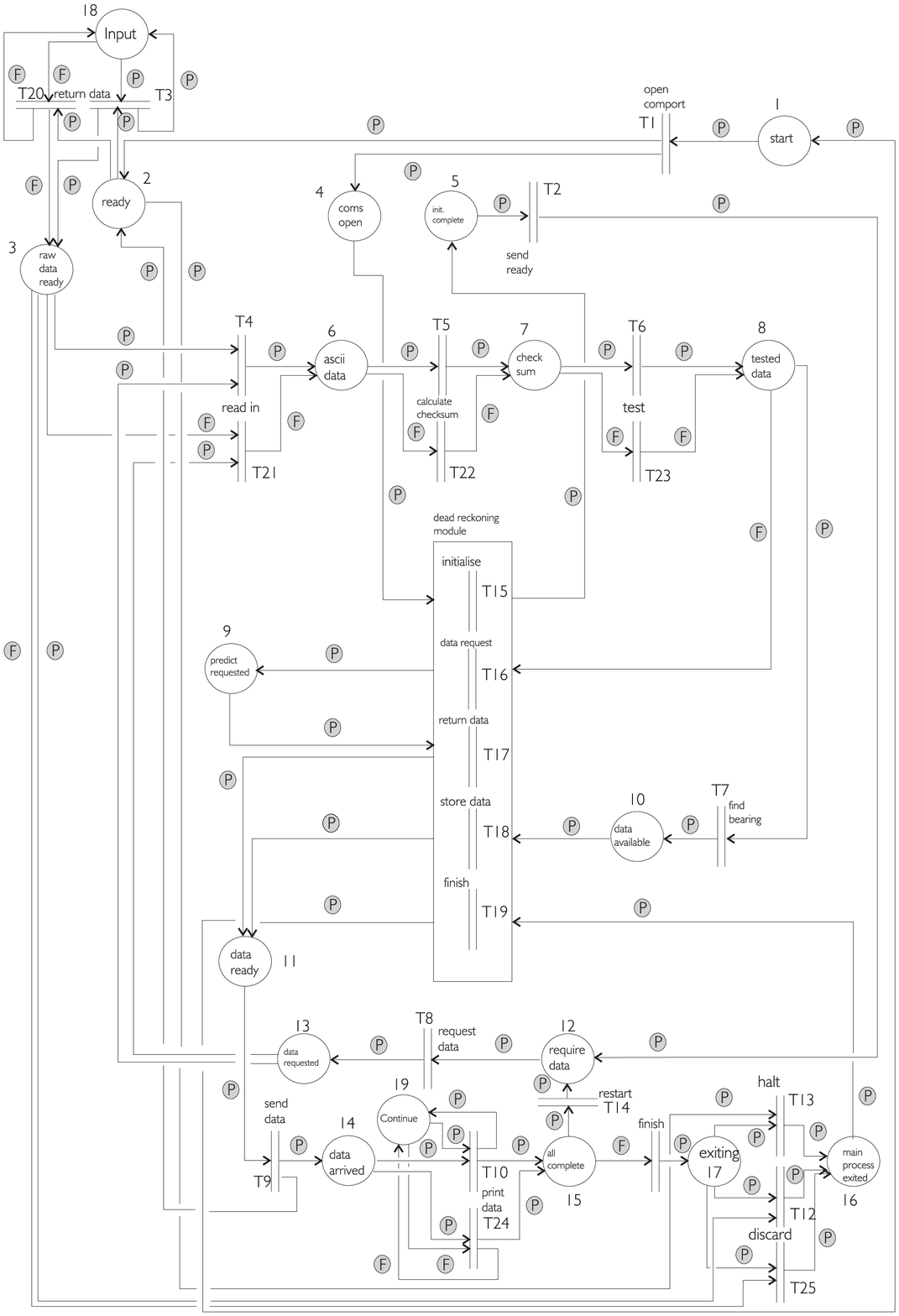}

\end{center}

\end{example}

\newpage
\section{Further Considerations}

\subsection{Boundedness}
Another interesting property is {\em boundedness} -- the maximum number 
of tokens that may exist at a particular place or 
the maximum number of tokens that can exist in the entire net --given an
initial marking. It is obvious to see how the catalogue may be used to check
either type of boundedness, but more interesting to observe that certain
information may be derived directly from the ($tdeg$) Gr\"obner basis.
If the Gr\"obner basis contains only polynomials $l-r$ (assume $l>r$) 
such that $l$ and $r$ are power products of the same total degree then all 
reachable markings will have the same number of tokens. 
The least number of tokens possible is the degree of the reduced form of
the polynomial associated with the initial marking.
Regarding the polynomials $l-r$ as reduction rules $l \to r$
we can sometimes determine the most number of tokens possible by examining the
degree-reducing rules to find what multiples of the reductum can be reduced
to the same form as the initial marking (it was possible to do this with the
47 rule Gr\"obner basis obtained for our last example). 
 
\subsection{Use and Efficiency}

Similarly to \cite{Caprotti} we point out that although in general Gr\"obner
basis computation can be lengthy, the type arising from Petri nets are not 
usually complex, involving only two-term polynomials with unitary coefficients.
There is no problem, in any case with ordinary or coloured Petri nets, as 
commutative Gr\"obner bases can always be found, using a computer algebra package
(e.g. $\mathsf{MAPLE}$).\\

Although it is possible to make use of existing implementations of
Buchberger's Algorithm it would be practical to include the Gr\"obner basis 
procedures as part of the software in our mechatronic navigation systems.
One aim of the research in \cite{Angie} is to provide an easier way of safely
programming a mobile robot. By using a Petri net to model the 
navigation system the C code controlling the robot is split into small pieces,
corresponding to the transitions in the net. A transition can be programmed in a
few lines, and code for a selection of alternative transitions could be provided
in advance. The structure of the net corresponds to the structure of the 
executable program, and thus by replacing individual transitions in the net 
the whole program for controlling the mobile robot can be rewritten and 
retested with the minimum difficulty. The Gr\"obner basis tests would form an
important part of the software, particularly in terms of safety.
One example this work could be applied to would be an autonomous
excavator. By using the Petri net representation, modifications to the control
of the excavator could be made in the field, without the requirement for on site
programming expertise. The Gr\"obner basis testing would provide a catalogue of 
reachable markings. If any undesirable (dangerous) states of the Petri net
were shown to be reachable, this problem could be rectified by further 
alteration to the net until the model was shown to be satisfactory.

\subsection{Streamed Petri Nets}

We are interested in Petri nets that can model systems involving
streams of data. 
Places will hold ordered lists of coloured tokens
rather than unordered sets of tokens. This introduces a degree of
noncommutativity into the Petri net. 
The Gr\"obner basis 
situation is more interesting here than with the ordinary Petri nets. 
Undecidability of the word problem \cite{TMora} indicates
the existence of streamed Petri nets for which it is not possible to determine 
whether or not a state is reachable.
The streamed models we have worked with store the streams of data as stacks
or allow random access to any substream of data within a given stream.
The problem with this is that the type of streamed Petri net suitable 
for our more advanced models is one 
whose transitions read data streams from the left and build them up on the 
right. This is a net to which we cannot yet apply Gr\"obner basis theory, 
but hope to investigate in future work.

\subsection{Enhanced Petri Nets}
Inhibitor arcs are the simplest extension to a basic Petri net.
The inhibitor arc is represented by a line with a small
circle at the end, equivalent to the {\small{$\mathtt{NOT}$}} in switching theory,
and is used to prevent a transition from firing.
If a transition $t$ has an inhibitor arc from a place $p$ then
$t$ is enabled only when there are tokens in all of its ordinary
input places and no tokens in the place $p$.
The inhibitor arcs provide an alternative method of forcing a decision
between two enabled transitions. These decisions can also be made
randomly, or with the use of colours, but in this specific case, the
inhibitor arc can give one transition priority over the other by
preventing the second transition from firing.  This method of decision
making could be useful in any system where one function should be given
priority over another.  For instance, if a Petri net driving a mobile
robot detected an obstruction, it would be important that it should
stop, or alter the speed of the motors before attempting to read any
sensors.\\

It is interesting to consider how the Gr\"obner
basis methods could be extended to cover variations of the Petri net theory, 
especially when the results of the extensions are motivated by the requirement
for testing modifications to navigation systems.

\subsection{Linked Petri Nets}
The motivation for our work has been the application to control systems of
mobile robots, using the TRAMP philosophy (Toolkit for Rapid Autonomous Mobile 
Prototyping). It allows the analysis of different control components of
a single mobile robot and it would be desirable for the Petri nets to be
logically linked to provide a unified model of the control of the device.
The analysis of the nets by Gr\"obner bases should then be extended to provide
an analysis of the model as a whole. The problem of the subdivision of
a large net into suitable components (objects) and the extension of local
analyses of such components to global checks on reachability, safety etc, are
examples of the well known local to global problem.

\end{document}